\documentclass[journal,usenames,dvipsnames]{IEEEtran}

\usepackage{myStyleIEEE}
\usepackage{bm}
\usepackage{smartdiagram}
\usesmartdiagramlibrary{additions}
\begin{document}
	
\title{Data-driven Local Control Design for Active Distribution Grids using off-line Optimal Power Flow and Machine Learning Techniques}

\author{
\IEEEauthorblockN{Stavros~Karagiannopoulos,~\IEEEmembership{Student Member,~IEEE,}
  Petros~Aristidou,~\IEEEmembership{Member,~IEEE,}\\
  and~Gabriela~Hug,~\IEEEmembership{Senior Member,~IEEE}}%
  \thanks{S. Karagiannopoulos and G. Hug are with the Power Systems Laboratory, ETH Zurich, 8092 Zurich, Switzerland. Email: \{karagiannopoulos $|$ hug\}@eeh.ee.ethz.ch.}%
  \thanks{P. Aristidou is with the School of Electronic and Electrical Engineering, University of Leeds, Leeds LS2 9JT, UK. Email: p.aristidou@leeds.ac.uk}%
}

\maketitle

\IEEEpeerreviewmaketitle
\begin{abstract}
The optimal control of distribution networks often requires monitoring and communication infrastructure, either centralized or distributed. However, most of the current distribution systems lack this kind of infrastructure and rely on sub-optimal, fit-and-forget, local controls to ensure the security of the network. In this paper, we propose a data-driven algorithm that uses historical data, advanced optimization techniques, and machine learning methods, to design local controls that emulate the optimal behavior without the use of any communication. \review{We demonstrate the performance of the optimized local control on a three-phase, unbalanced, low-voltage, distribution network. The results show that our data-driven methodology clearly outperforms standard industry local control and successfully imitates an optimal-power-flow-based control.} 
\end{abstract}


\section{Introduction}
Some of the most notable developments foreseen in power systems target Distribution Networks (DNs). In the future, DNs will host a large percentage of Distributed Generators (DGs), including Renewable Energy Sources (RES), to supply a growing share of the total demand. These units, in combination with other Distributed Energy Resources (DERs) such as electric vehicles, Battery Energy Storage Systems (BESSs) and Flexible Loads (FLs), will elevate the role of Distribution System Operators (DSOs), allowing them to provide ancillary services and support the bulk transmission system~\cite{PES-TR22}. However, this new paradigm introduces significant challenges to the DN operation~\cite{PES-TR22}.

Traditionally, to address these challenges, DSOs have relied only on grid reinforcement and ignored the flexibility offered by DERs. This approach is now unable to cope with the new challenges while keeping the cost for the consumer low and achieving high security and reliability goals. It is apparent  that DSOs need to operate DNs actively, involving DERs to ensure secure, reliable and cost-effective operation.

Based on the communication infrastructure available for controlling the DERs, operational schemes can be broadly classified as centralized, distributed and decentralized or local. Centralized schemes require extensive monitoring and communication infrastructure and usually leverage the performance of powerful optimization-based control techniques. The capabilities of extensive monitoring and communication infrastructure allow for system-wide optimal operation by coordinated control of DERs~\cite{Fortenbacher2016a,KaragiannopoulosGM}.  This type of control has lately attracted significant attention thanks to advances in computational power, wireless communication, and new theoretical developments in approximations of the nonlinear AC power flow equations~\cite{Lavaei2012,Molzahn2015}.  A lot of methods rely on semi-definite relaxations, e.g.~\cite{Lavaei2012}, which find global optimal solutions in many practical cases with specific conditions, but not in the general case~\cite{Molzahn2014}. Lately, many researchers started dealing with multi-phase systems, e.g.~\cite{Zamzam2018,Bolognani_manifoldlinearization}. A very efficient method is presented in~\cite{Bolognani_manifoldlinearization} based on linear manifold approximants, while in~\cite{Zamzam2018} the authors use an iterative algorithm to solve the OPF as a nonconvex quadratically constrained quadratic program. However, they do not model explicitly the power losses in LV grids, neither consider uncertainties.

Nevertheless, the infrastructure required for this type of control is rarely available in DNs, and the financial benefit for investing in such capabilities not clear. Decentralized control strategies, e.g.~\cite{Tonkoski2011,Kotsampopoulos2013}, tackle power quality and security problems using only local measurements. These type of controls are widely used in DNs today and have been embedded in several grid codes. The benefit of these methods lies in the simplicity and the relatively low cost of implementation. No communication infrastructure is needed, keeping the required investment at a minimum. However, these methods usually employ a one-size-fits-all approach, where the same control parameters are employed in all DNs, different generator types, and operating conditions. This approach can lead to unforeseen problems, especially in a rapidly changing environment.

Finally, distributed approaches, e.g.~\cite{Bolognani2013,OAEC2016ja}, use limited communication between different DERs to coordinate them and achieve a close-to-optimal operation. While these methods try to bridge the gap between local and centralized methods, they still require some communication infrastructure and usually employ consensus-based control algorithms which are sensitive to communication delays and errors.

\begin{figure*}
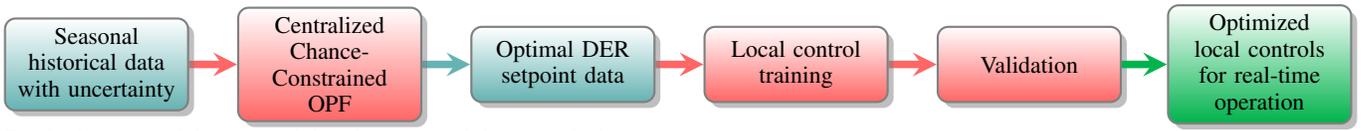

\begin{centering}
\smartdiagramset{%
  set color list={teal!60,red!60,teal!60,red!60,red!60},
  module x sep=3.1,
  text width=2.2cm,
  back arrow disabled=true,
  additions={
    additional item offset=0.6cm,
    additional item border color=gray,
    additional item font=\small,
    additional arrow color=teal!60!green,
    additional item text width=2.2cm,
    additional item bottom color=teal!60!green,
    additional item shadow=drop shadow,
  }
}
\smartdiagramadd[flow diagram:horizontal]{%
  Seasonal historical data with uncertainty, Centralized Chance-Constrained OPF, Optimal DER setpoint data, Local control training,Validation}{right of module5/Optimized local controls for real-time operation}
\smartdiagramconnect{->}{module5/additional-module1}
\end{centering}
\caption{Overview of the proposed data-driven control design method.}\label{fig:overview1} 
\end{figure*}

Lately, data-driven methods have attracted a lot of attention in the power systems area~\cite{Guo2018a,Guo2018b,Garg2018,Dobbe2018,Sondermeijer16}.
In~\cite{Guo2018a}, data-based methods are used to solve a distributionally robust OPF problem. The method is based on a model predictive control algorithm that utilizes forecast error training datasets, and the focus is on obtaining closed-loop control policies which are robust to sampling errors in the data. In~\cite{Guo2018b}, the authors demonstrate the method to mitigate overvoltages, assuming balanced phase loading. However, this approach requires a centralized scheme as well as reliable and accurate monitoring and communication infrastructure. \review{Using machine learning techniques to represent the optimal behavior is studied  in~\cite{Garg2018,KaragiannopoulosIET,FedericaGM,Dobbe2018,Sondermeijer16}.
Reference~\cite{Garg2018} uses non-linear control policies to calculate the real-time reactive power injections of the inverter-based DGs.} Although it uses a linearized version of the grid, assumes balanced operation and focuses only on one measure, i.e. reactive power control, this scheme is very flexible due to the various kernel functions which are able to model complex and non-linear behaviors.
In~\cite{Sondermeijer16} and~\cite{Dobbe2018}, multiple linear regression is used in an open-loop fashion to calculate a function for each inverter that maps its local historical data to pre-calculated optimal reactive power injections. However, both references consider only reactive power control, neglecting possible combinations with other available controls, and reference~\cite{Sondermeijer16} assumes a balanced DN, i.e. using a single-phase representation.

The focus of this work is on distribution grids where communication and monitoring infrastructures have not been deployed yet. Thus, in this paper, we propose a data-driven control design method to derive optimized \textit{local} controls for several types of DERs. However, it should be mentioned that even in an environment with deep penetration of communication capabilities, local controls are valuable as they provide a scalable approach to efficiently utilize an increasingly large number of distributed resources. The methodology is sketched in Fig.~\ref{fig:overview1} and detailed hereafter. 

First, we use a model of the DN under study along with historical generation and consumption data. Then, we employ an offline centralized optimization algorithm to compute the \textit{optimal} DER control setpoints for different operating conditions. The objective of the offline algorithm is to minimize the system losses and adjustments of DER resources while ensuring system security and power quality. The formulation takes into account the uncertainty coming from RES and the unbalanced, three-phase, operation. Finally, we use Machine Learning (ML) techniques applied on the optimal setpoints obtained from this optimization to design local DER controls for the real-time operation of the DN. In this way, we derive simple and efficient optimized local controls that can mimic the behaviour of centralized optimization-based schemes, without the need of any communication infrastructure.

This paper extends and completes our previous work in~\cite{KaragiannopoulosIET,stavrosPowertech,FedericaGM}. In~\cite{KaragiannopoulosIET}, we presented the idea of designing customized control schemes for each DER based on off-line centralized optimization. However, we only considered reactive power control and active power curtailment, and the derivation of the local volt-var curves was not based on ML techniques, but on rule-based heuristics. In~\cite{stavrosPowertech}, we added controllable loads in our methodology with a simple rule-based real-time control scheme. \review{Finally, in~\cite{FedericaGM} we utilized Support Vector Machines (SVMs) to derive local volt-var curves, but we accounted ex-post for the needed monotonicity and slope constraints of the final curves.}

In this paper, we consider reactive power control, active power curtailment, controllable load shifting and battery energy storage systems. The local schemes are derived by several machine learning techniques, such as segmented regression and SVMs as regressors and classifiers.  
More specifically, the contributions of this paper can be summarized as follows:
\begin{itemize}
    \item A computationally tractable off-line centralized control algorithm based on a three-phase, multi-period, Chance-Constrained Optimal Power Flow (CC-OPF), considering RES uncertainty and unbalanced operation.
    \item A novel data-driven local control design methodology for the optimal operation of several types of DERs, using different regression and classification ML techniques.
\end{itemize}

It should be noted that in this work, we use the centralized scheme \emph{off-line}, i.e. it does not require real-time monitoring and communication infrastructure. Instead, we use historical values collected in the past. Although any OPF formulation suitable for DNs can be used to derive the optimal DER setpoints, the proposed formulation allows us to use a \emph{tractable} three-phase multi-time OPF formulation that can consider \emph{uncertainties} and various models of DERs.


The remainder of the paper is organized as follows: In Section~\ref{OPFall}, we present the mathematical formulation of the CC-OPF algorithm used to obtain the optimal DER setpoints. Then, in Section~\ref{decentralized}, we describe the ML methods used for deriving the optimized local control schemes of the DERs. In Section~\ref{case}, we introduce the case study and simulation results that show the performance of the optimized controllers. Finally, we draw conclusions in Section~\ref{Conclusion}.

\section{Centralized Chance-Constrained OPF}\label{OPFall}

In this section, the centralized CC-OPF scheme used to compute the optimal DER setpoints for different operating conditions is presented. The objectives and constraints of the OPF-based algorithm are vital for the overall methodology as they will be reflected in the generated optimal DER setpoint data and will in turn influence the local control design. 

It should be noted that having enough data to run an OPF is critical to the process. The data can be gathered using low-cost energy monitoring devices and if some data are missing or noisy, we can extrapolate using historical data, public domain information, or information from neighbouring systems.

\vspace{-0.1cm}
\subsection{Centralized OPF}\label{detOPF}

\subsubsection{Objective function}\label{objfun}

The objective function selected includes minimizing the cost of DER control and the network losses, over all of the network nodes ($N_b$), phases ($z$) and branches ($N_{br}$) for the entire time horizon ($N_{hor}$). This is described by: 
\begin{align}
    \min_{ \bm{u} } & \sum \limits_{t=1}^{N_\textrm{hor}} \biggl\{ \sum \limits_{z\in\{a,b,c\}} \sum \limits_{j=1}^{N_\textrm{b}}  \biggl(C_\textrm{P} \hspace{-0.05cm} \cdot \hspace{-0.05cm} P_\textrm{curt,j,z,t} \hspace{-0.05cm} + \hspace{-0.05cm} C_\textrm{Q} \hspace{-0.05cm} \cdot \hspace{-0.05cm} Q_\textrm{ctrl,j,z,t} \biggr) \hspace{-0.05cm} \nonumber  \\  
    & + \sum \limits_{z\in\{a,b,c\}} \sum \limits_{i=1}^{N_\textrm{br}} C_\textrm{P} \hspace{-0.05cm} \cdot \hspace{-0.05cm} P_\textrm{loss,i,z,t} \biggl\} \cdot \Delta t \label{eq:objfun1} \\
    &+ C_\textrm{H} \cdot \biggl( ||\eta_\textrm{V}||_{\infty} + ||\eta_{\textrm{I}}||_{\infty} +||\eta_{\textrm{VUF}}||_{\infty} \biggr) \nonumber 
\end{align}
\noindent where $\bm{u}$ is the vector of the available active control measures and $\Delta t$ is the length of each time period. The curtailed power of the DGs connected at phase $z$, at node $j$ and time $t$ is given by $P_\textrm{curt,j,z,t} = P_\textrm{g,j,z,t}^{\textrm{max}} - P_\textrm{g,j,z,t}^{\textrm{ }}$, where $P_\textrm{g,j,z,t}^{\textrm{max}}$ is the maximum available active power and $P_{\textrm{g,j,z,t}}^{\textrm{ }}$ the active power injection of the DGs. The use of reactive power support $Q_\textrm{ctrl,j,z,t}= | Q_\textrm{g,j,z,t}^{\textrm{ }}|$ for each DG connected to phase $z$ of node $j$ and time $t$ is also minimized; $Q_{\textrm{g,j,z,t}}^{\textrm{ }}$ represents the DG reactive power injection or absorption. The coefficients $C_\textrm{P}$ and $C_\textrm{Q}$ represent, respectively, the DG cost of curtailing active power and providing reactive power support (DG opportunity cost or contractual agreement). The assumption that $C_\textrm{Q} \ll C_\textrm{P}$ is made, which prioritizes the use of reactive power control over active power curtailment. 
In our case, we follow the method of~\cite{Strunz2014} and perform Kron's reduction in order to use a three-phase three-wire power flow representation. In this case, the phase voltages and currents are obtained with acceptable accuracy, but as it is shown in~\cite{Ochoa2005}, calculating the losses using the current magnitude squared times the resistance formula, leads to high overestimation of the correct total losses. Thus, in this case, one can calculate the total losses, by using the difference between input and output power in each phase \cite{Kersting2002}. \review{Thus, $P_\textrm{loss,i,z,t} =  |\textrm{Re} (S_{\textrm{$\textrm{i}_f$,z,t}} + S_{\textrm{$\textrm{i}_t$,z,t}})|$, where $S_{\textrm{$\textrm{i}_f$,z,t}}$ and $S_{\textrm{$\textrm{i}_t$,z,t}}$ represent the apparent power flowing into branch $i$ from each end; $\textrm{i}_f$ and $\textrm{i}_t$ are the sending and receiving ends of the branch.} Finally, $C_\textrm{H}$ is a large penalty associated with violating the security and power quality constraints. It is used in conjunction with the variables $(\eta_\textrm{V},\eta_{\textrm{I}},\eta_{\textrm{VUF}})$ to relax respectively the voltage, thermal or balancing constraints and avoid infeasibility. When one of these limits is binding, the output of the overall objective function gets dominated by this term and might lose a real monetary meaning (unless the cost of violating the security and power quality constraints is quantified and monetized by the DSO).

\subsubsection{Power balance constraints}\label{powerbalance}
The power injections at every node $j$, phase $z$ and time step $t$ are given by
\begin{subequations} \label{eq:node_balance}
\begin{align}
	P_{\textrm{inj,j,z,t}}^{\textrm{ }}&=P_{\textrm{g,j,z,t}}^{\textrm{ }} - P_{\textrm{lflex,j,z,t}}^{\textrm{ }} - (P_{\textrm{B,j,z,t}}^{\textrm{ch}} - P_{\textrm{B,j,z,t}}^{\textrm{dis}}) \label{eq:node_balance_P}\\
    Q_{\textrm{inj,j,z,t}}^{\textrm{ }}&=Q_{\textrm{g,j,z,t}}^{\textrm{ }} - P_{\textrm{lflex,j,z,t}}^{\textrm{ }} \cdot \textrm{tan}(\phi_{\textrm{load}}) + Q_{\textrm{B,j,z,t}}  \label{eq:node_balance_Q}   
\end{align}
\end{subequations}
where $P_{\textrm{lflex,j,z,t}}^{\textrm{ }}$ and $P_{\textrm{lflex,j,z,t}}^{\textrm{ }} \cdot \textrm{tan}(\phi_{\textrm{load}})$ are the active and reactive node demands (after control) of constant power type, with $\textrm{cos}(\phi_{\textrm{load}})$ being the power factor of the load; $Q_{\textrm{B,j,z,t}}$ the reactive power of the BESS and, $P_{\textrm{B,j,z,t}}^{\textrm{ch}}$ and $P_{\textrm{B,j,z,t}}^{\textrm{dis}}$ are respectively the charging and discharging BESS active powers.

\subsubsection{Power flow constraints}\label{powerflow}
The non-linear AC power-flow equations that model the DN network make solving the OPF problem computationally challenging. Since the OPF will be used to process several scenarios in a multi-period framework, it is necessary to use some approximations to increase its computational performance. For this reason, the iterative Backward/Forward Sweep (BFS) power flow~\cite{Teng2003} method is used in this work, extending the formulation presented by the authors in \cite{stavrosPowertech,StavrosIREP,StavrosPSCC18} for a three-phase, unbalanced system.

Following our previous work~\cite{StavrosPSCC18}, a single iteration of the BFS power-flow method is used to replace the AC power-flow constraints in the OPF formulation. This is written as ($j=1,\ldots N_{b}$, $z \in \{a,b,c\}$):
\begin{gather}
    I_\textrm{inj,j,z,t}= \left (\frac{(P_{\textrm{inj,j,z,t}}^{\textrm{ }} + jQ_{\textrm{inj,j,z,t}}^{\textrm{ }})^{*}}{\bar{V}_{\textrm{j,z,t}}^{*}}\right) \nonumber  \\
    \bm{I}_\textrm{br,t}=\bm{BIBC} \cdot \bm{I}_\textrm{inj,t} \nonumber \\
    \Delta{\bm{V}}_{\textrm{t}}=\bm{BCBV} \cdot \bm{I}_\textrm{br,t} \nonumber  \\
    \bm{V}_{{\textrm{t}}}= \bm{V}_\textrm{slack} - \Delta V_\textrm{tap} \cdot \rho_{\textrm{t}} + \Delta{\bm{V}}_{\textrm{t}} \nonumber \\
    \rho_{min} \leq \rho_{\textrm{t}} \leq \rho_{max}, \label{eq:OLTC3}  
\end{gather}
\noindent where $\bar{V}_{\textrm{j,z,t}}^{*}$ is the voltage of phase $z$, at node $j$ at time $t$, $~^*$~indicates the complex conjugate and the bar indicates that the value from the previous iteration is used (details will be given later); $\bm{I}_\textrm{inj,t}^{\textrm{}}$ and $\bm{I}_\textrm{br,t}^{\textrm{}}$ are respectively the vectors of the three-phase bus injection and branch flow currents; and, $\bm{BIBC}$ (Bus Injection to Branch Current) is a matrix with ones and zeros, capturing the three-phase topology of the DN (including any single-phase laterals); $\Delta{\bm{V}}^{}_{\textrm{t}}$ is the vector of voltage drops over all branches and phases; $\bm{BCBV}$ (Branch Current to Bus Voltage) is a matrix with the complex impedance of the lines as elements (including mutual coupling); $\bm{V}_\textrm{slack}$ is the three-phase voltage in per unit at the slack bus (here assumed to be $\{1\hspace{-0.1cm}<\hspace{-0.1cm}0^{\circ{}}, 1\hspace{-0.1cm}<\hspace{-0.1cm}-120^{\circ{}},1\hspace{-0.1cm}<\hspace{-0.1cm}120^{\circ{}}\}$); $\Delta V_{tap}$ is the voltage magnitude change caused by one tap action of the On-Load Tap Changer (OLTC) transformer and assumed constant for all taps for simplicity; and, $\rho_{\textrm{t}}$ is an integer value defining the position of the OLTC position. The parameters ($\rho_\textrm{min},\rho_\textrm{max}$) are respectively the minimum and maximum tap positions of the OLTC transformer.

This convex formulation provides a good approximation to the nonlinear AC OPF~\cite{Fortenbacher2016a}, is computationally tractable even in a three-phase model~\cite{StavrosPSCC18}, and results in AC feasible solutions which can account for uncertainties, see~\cite{StavrosIREP} and Section~\ref{uncertain}.

\subsubsection{Thermal loading and voltage constraints} \label{VIconstraints}
The constraint for the current magnitude for  branch $i$ and phase $z$ at time $t$ is given by
\begin{align}
    |I_\textrm{br,i,z,t}|   \leq I_{\textrm{i,z,max}} + \eta_{\textrm{I,i,z,t}} \label{eq:Qk1}, \qquad 
    \eta_{\textrm{I,i,z,t}}  \geq 0
\end{align}
\noindent where $I_\textrm{br,i,z,t}$ is the branch current; $I_{\textrm{i,z,max}}$ is the maximum thermal limit; and, $\eta_{\textrm{I,i,z,t}}$ is used to relax the constraint when the thermal constraints cannot be met.

Similarly, the voltage constraints are given by
\begin{align}
    V_\textrm{min} - \eta_{\textrm{V,j,z,t}} & \leq | V_{\textrm{j,z,t}} | \leq V_{\textrm{max}} +\eta_{\textrm{V,j,z,t}} \label{eq:V1}, \qquad
    \eta_{\textrm{V,j,z,t}} & \geq 0
\end{align}
\noindent where $(V_{\textrm{max}}$, $V_{\textrm{min}})$ are respectively the upper and lower acceptable voltage limits and $\eta_{\textrm{V,j,z,t}}$ is used to relax the constraint when the voltage constraints cannot be met.

Unfortunately, \eqref{eq:V1} is non-convex due to the minimum voltage magnitude requirement. In order to avoid the non-convexity, we rotate the three voltage phases $\{a, b, c\}$ by~{$\mathcal{R} = \{1\hspace{-0.1cm}<\hspace{-0.1cm}0^{\circ{}}, 1\hspace{-0.1cm}<\hspace{-0.1cm}120^{\circ{}},1\hspace{-0.1cm}<\hspace{-0.1cm}-120^{\circ{}}\}$} so that they lie close to the reference axis $0^{\circ{}}$ and we define the same feasible space for each of the three phases (see~\cite{StavrosPSCC18} for more details)
\begin{align}
    \begin{cases}
               |\mathcal{R} V_{\textrm{j,z,t}} | \leq V_{\textrm{max}} + \eta_{\textrm{V,j,z,t}} \\
               \textrm{Re}\left\{\mathcal{R} V_{\textrm{j,z,t}}    \right\} \geq V_\textrm{min} - \eta_{\textrm{V,j,z,t}} 
            \end{cases}\label{eq:V4}
\end{align}

\subsubsection{Balancing constraint}\label{Unbalance}
A balancing constraint is used to improve the power quality of the DN by balancing the three phase voltages. We use the IEC unbalance definition~\cite{Pillay2001,IEC} of Voltage Unbalance Factor ($VUF$), which is given by $VUF(\%)=100\%\frac{|V_{-}|}{|V_{+}|}$, where $V_{-}$ and $V_{+}$ are respectively the negative and positive sequence derived by symmetrical component analysis. 

The balancing constraint for node $j$ and time $t$ is given by $VUF_{\textrm{j,t}}(\%) \leq VUF_{\textrm{MAX}}$, where $VUF_{\textrm{MAX}}$ is the acceptable voltage unbalance factor (e.g. $2\%$ for $95\%$ of the week according to EN50160~\cite{EN50160}). Since this constraint is non-convex, we approximate $VUF$ by the negative voltage sequence~\cite{StavrosPSCC18}, assuming the positive voltage sequence is very close to $1$~pu, i.e. $VUF_{\textrm{j,t}} \approx |V_{\textrm{-,j,t}}|$. This gives
\label{eq:VUFall}
\begin{align}
    |V_{\textrm{-,j,t}}| \leq VUF_{\textrm{MAX}} + \eta_{\textrm{VUF,j,t}}\label{eq:VUF} ,\ \qquad   \eta_{\textrm{VUF,j,t}}  \geq 0
\end{align}
\noindent where $\eta_{\textrm{VUF,j,t}}$ relaxes the constraint when it cannot be met.

\subsubsection{DER constraints}\label{DERineq}

\paragraph{DG limits}
In this work, without loss of generality, we only consider inverter-based DGs such as PVs. Their limits are thus given by
\begin{subequations}
\label{eq:prod3}
\begin{align}
    P_{\textrm{g,j,z,t}}^{\textrm{min}} \leq P_{\textrm{g,j,z,t}}^{\textrm{ }} \leq P_{\textrm{g,j,z,t}}^{\textrm{max}}\\ Q_{\textrm{g,j,z,t}}^{\textrm{min}} \leq  Q_{\textrm{g,j,z,t}}^{\textrm{ }} \leq Q_{\textrm{g,j,z,t}}^{\textrm{max}}
\label{eq:PV_prod3Qa}
\end{align}
\end{subequations}
\noindent where $P_{\textrm{g,j,z,t}}^{\textrm{min}}$, $P_{\textrm{g,j,z,t}}^{\textrm{max}}$, $Q_{\textrm{g,j,z,t}}^{\textrm{min}}$ and $Q_{\textrm{g,j,z,t}}^{\textrm{max}}$ are the upper and lower limits for active and reactive DG power at each node $j$, phase $z$ and time $t$. These limits vary depending on the type of the DG and the control schemes implemented. Usually, small DGs have technical or regulatory~\cite{VDE} limitations on the power factor they can operate at or reactive power they can produce. This restriction can be included by linking the active and reactive power limits in \eqref{eq:prod3} through the maximum power factor value. 

\paragraph{Controllable loads} 
Moreover, we consider flexible loads which can shift a fixed amount of energy consumption in time. The behavior of the loads is given by
\begin{equation}\label{eq:CL}
    P_{\textrm{lflex,j,z,t}}^{\textrm{ }} = P_{\textrm{l,j,z,t}}^{\textrm{ }} + n_{\textrm{j,z,t}} \cdot P_{\textrm{shift,j,z}}, \qquad
   \sum \limits_{t=1}^{N_{hor}} n_{\textrm{j,z,t}}=0
\end{equation}
\noindent where $P_{\textrm{lflex,j,z,t}}^{\textrm{ }}$ is the controlled active power demand at phase $z$ of node $j$ and at time $t$, $P_{\textrm{shift,j,z}}$ is the load that can be shifted (assumed constant) and $n_{\textrm{j,z,t}} \in \left\{-1,0,1\right\}$ is an integer variable indicating an increase or a decrease of the load when shifted from the initial demand $P_{\textrm{l,j,z,t}}^{\textrm{ }}$. We assume that the final total daily energy demand needs to be maintained.

\paragraph{Battery Energy Storage Systems}
Finally, the constraints related to the BESS are given as
\begin{subequations}
\label{eq:BESS}
\begin{gather}
    SoC_{\textrm{min}}^{\textrm{bat}} \cdot E_{\textrm{cap,j,z}}^{\textrm{bat}} \leq  E_{\textrm{j,z,t}}^{\textrm{bat}} \leq SoC_{\textrm{max}}^{\textrm{bat}} \cdot E_{\textrm{cap,j,z}}^{\textrm{bat}}      \label{eq:BESS_en}\\  
    E_{\textrm{j,z,1}}^{\textrm{bat}} = E_{\textrm{start}}    \label{eq:BESS_SOC}\\
    E_{\textrm{j,z,t}}^{\textrm{bat}}  = E_{\textrm{j,z,t-1}}^{\textrm{bat}} + (\eta_{\textrm{bat}} \cdot  P_{\textrm{B,j,z,t}}^{\textrm{ch}} - \frac{P_{\textrm{B,j,z,t}}^{\textrm{dis}}}{\eta_{\textrm{bat}}}) \cdot \Delta t    \label{eq:dynBESS1}\\
    0 \leq P_{\textrm{B,j,z,t}}^{\textrm{ch}} \leq P_{\textrm{max}}^{\textrm{bat}} ,  \quad  0 \leq P_{\textrm{B,j,z,t}}^{\textrm{dis}} \leq P_{\textrm{max}}^{\textrm{bat}}  \label{eq:BESS_P} \\
    P_{\textrm{B,j,z,t}}^{\textrm{ch}} +   P_{\textrm{B,j,z,t}}^{\textrm{dis}} \leq  \textrm{max}(P_{\textrm{B,j,z,t}}^{\textrm{ch}}, P_{\textrm{B,j,z,t}}^{\textrm{dis}}) \label{eq:BESS_P2} \\
    Q_{\textrm{B,j,z,t}}^{2} \leq  (S_{\textrm{max}}^{\textrm{bat}})^{2} - \textrm{max}((P_{\textrm{B,j,z,t}}^{\textrm{ch}})^{2}, (P_{\textrm{B,j,z,t}}^{\textrm{dis}})^{2})
\end{gather}
\end{subequations}
where $E_{\textrm{cap,j,z}}^{\textrm{bat}}$ is the installed BESS capacity connected to phase $z$ at node $j$; 
$SoC_{\textrm{min}}^{\textrm{bat}}$ and $SoC_{\textrm{max}}^{\textrm{bat}}$ are the fixed minimum and maximum per unit limits for the battery state of charge; and, $E_{\textrm{j,z,t}}^{\textrm{bat}}$ is the available energy at node $j$, phase $z$ and time $t$. The initial energy content of the BESS in the first time period is given by $E_{\textrm{start}}$, and (\ref{eq:dynBESS1}) updates the energy in the storage at each period $t$ based on the BESS efficiency $\eta_{bat}$, time interval $\Delta t$ and the charging and discharging power of the BESS $P_{\textrm{B,j,z,t}}^{\textrm{ch}}$ and $ P_{\textrm{B,j,z,t}}^{\textrm{dis}}$. The charging and discharging powers are defined as positive according to \eqref{eq:BESS_P}, while \eqref{eq:BESS_P2} is re-casted as mixed-integer constraint with 2 binaries for each time step, and  ensures that the BESS is not charging and discharging at the same time.

\subsection{Accounting for Uncertainty through Chance Constraints}\label{uncertain}
To account for the effect of generation uncertainty and to limit possible adverse effects on the security constraints, we reformulate the problem using chance constraints. Chance-constrained optimization problems aim to keep the probability of certain random events below targeted values. In this work, we assume that the PV power injection is the only source of uncertainty. However, load uncertainty can be also included in a similar way. The interested reader is referred to~\cite{Bienstock2014} for a general overview of risk-aware control under uncertainties.

\subsubsection{Formulation of the Chance Constraints}
The branch current flows and the voltage magnitudes are functions of the power injections and are hence directly influenced by the PV power uncertainty. Thus, we model the corresponding voltage and current constraints as chance constraints that will hold with a chosen probability $1-\varepsilon$, where $\varepsilon$ is the acceptable violation probability. E.g., the maximum voltage magnitude constraint is reformulated as $\mathbb{P}\left\{|V_{\textrm{bus,j,t}}| \leq  V_{\textrm{max}} \right\}\geq{1-\varepsilon}$\cite{Roaldirep}.

To solve the resulting CC-OPF, we need to reformulate the constraints in a tractable form. This can be achieved by using an analytical form assuming a certain distribution of the forecast error~\cite{roaldArxiv}, or a distribution-agnostic method~\cite{Guo2018a,roaldArxiv}. In this work, we follow~\cite{Roald2013,roald2017submitted} that rely on an iterative solution scheme that fits very well with the iterative nature of the BFS-OPF. The core idea is that the chance-constrained problem can be cast as a deterministic problem with tightened constraints. The tightenings represent security margins against uncertainty, i.e., uncertainty margins, which drive the trade-off between cost and system security, and are functions of the optimization variables. The iterative scheme alternates between solving the deterministic problem with a given set of tightenings, and evaluating the optimal constraint tightening based on the solution of the deterministic problem~\cite{StavrosIREP}. A feasible solution is found when the tightenings do not change between iterations.  Thus, we interpret the probabilistic constraints as tightened deterministic versions of the original constraints, and we express~\eqref{eq:Qk1} and~\eqref{eq:V4} \review{as
\begin{align}
        \begin{cases}
               |\mathcal{R} V_{\textrm{j,z,t}} | & \leq V_{\textrm{max}} - \Omega_\textrm{V j,z,t}^\textrm{upper} + \eta_{\textrm{V,j,z,t}} \\
               \textrm{Re}\left\{\mathcal{R} V_{\textrm{j,z,t}}    \right\} &\geq V_\textrm{min} + \Omega_\textrm{V j,z,t}^\textrm{lower} -\eta_{\textrm{V,j,z,t}} 
    \end{cases} \label{eq:vol_lim2}\\
    |I_\textrm{br,i,z,t}|  \leq I_{\textrm{i,z,max}} - \Omega_{I_\textrm{br,i}} + \eta_{\textrm{I,i,z,t}} \label{eq:current_lim2}
\end{align}
where $\Omega_\textrm{V}^\textrm{lower},~\Omega_\textrm{V}^\textrm{upper}$ are the tightenings for the lower and upper voltage magnitude constraints and $\Omega_{I_\textrm{br}}$ are the tightenings of the current magnitude constraints. The procedure is explained in more detail in~\cite{StavrosIREP}.}

\subsubsection{Uncertainty margin evaluation based on Monte Carlo Simulations}

to evaluate the uncertainty margins, we use a Monte Carlo method. The uncertainty margins are considered constant within the OPF solution, but then evaluated outside of the OPF iterations. The advantages of this method lie in the ability to use the non-linear AC power-flow and to have any uncertainty probability distribution. 

First, empirical distributions for the voltage and current chance constraints are formed at each time step based on the Monte Carlo simulations. To enforce a chance constraint with $1-\epsilon$ probability we need to ensure that the $1-\epsilon$ quantile of the distribution remains within the bounds. Thus, the tightening corresponds to the difference between the forecasted value with zero forecast error and the $1-\epsilon$ quantile value evaluated based on the empirical distribution resulting from the Monte Carlo Simulations, e.g. $|V_{\textrm{bus,j,t}}^{\textrm{0}}|$ and $|V_{\textrm{bus,j,t}}^{\textrm{1-$\epsilon$}}|$ for the voltage constraints.
The empirical uncertainty margins to be used in the next iteration are then given by 
\begin{subequations}
\begin{align}
    \Omega_\textrm{V j,t}^\textrm{upper}   &= |V_{\textrm{bus,j,t}}^{\textrm{1-$\epsilon$}}| - |V_{\textrm{bus,j,t}}^{\textrm{0}}|\\
    \Omega_\textrm{V j,t}^\textrm{lower}   &= |V_{\textrm{bus,j,t}}^{\textrm{0}}|  - |V_{\textrm{bus,j,t}}^{\textrm{$\epsilon$}}|\\
    \Omega_{I_\textrm{br,i}}^\textrm{upper}&= |I_\textrm{br,i,t}^{\textrm{1-$\epsilon$}}| - |I_\textrm{br,i,t}^{\textrm{0}}|
\end{align}
\end{subequations}
\review{where superscript $^0$ indicates the values at the operating point with zero forecast error.}

\subsubsection{Iterative Solution Algorithm}
Since the uncertainty margins rely on the selected DER setpoints, an iterative algorithm is used to solve the problem~\cite{Schmidli2016, roald2017submitted}. It alternates between solving a deterministic OPF with tightened constraints, and calculating the uncertainty margins $\Omega_\textrm{V}^\textrm{lower},~\Omega_\textrm{V}^\textrm{upper},~\Omega_{I_\textrm{br}}^\textrm{upper}$. When the change in the tightening values between two subsequent iterations is below a threshold $(\eta_V^{\Omega},~\eta_I^{\Omega})$, then the algorithm has converged.

\subsection{Solution Algorithm}
\label{SolAlgo}
\begin{figure}[]
\centering 
    \tikzstyle{note} = [rectangle, dashed, draw, fill=ProcessBlue, font=\footnotesize,
    text width=10cm, text centered, rounded corners, minimum height=6.7cm,minimum width=7.85cm,opacity=0.2]
\tikzstyle{note2} = [rectangle, dashed, draw, fill=YellowGreen, font=\footnotesize,
    text width=10cm, text centered, rounded corners, minimum height=13cm,minimum width=8.7cm,opacity=0.2]

\centering
\begin{tikzpicture}[node distance=2cm, align=center, scale=0.65, every node/.style={transform shape}]
        \node (start) [startstop, opacity=0.9] {{\large \textbf{Initialize:}} \\ $k=0$, $V_{\textrm{bus}}^{k}=\{1\hspace{-0.1cm}<\hspace{-0.1cm}0^{\circ{}}, 1\hspace{-0.1cm}<\hspace{-0.1cm}-120^{\circ{}},1\hspace{-0.1cm}<\hspace{-0.1cm}120^{\circ{}}\}$ \\
                                                            $m=1$, $\Omega_{i_{br}}^{m-1}=\Omega_{V i}^{m-1}=0$ };
        \node (in1) [process, below of=start, fill=Gray!30] {{\large Run three phase multi-period OPF\\ with one BFS iteration }};
        \node (pro1) [process, below of=in1 , fill=Gray!30] {{\large Run BFS power flow until convergence }};
        \node (dec1) [decision, below of=pro1, yshift=-0.5cm, , fill=Gray!30] {$max|(|V_{\textrm{bus}}^{k}| - |V_{\textrm{bus}}^{PF}|)|\leq \tilde{\eta}$};
        \node (pro2a) [process, below of=dec1, yshift=-0.5cm] {{\large Evaluate $\Omega_{V i}^{m}$, $\Omega_{I_{br}}^{m}$} \\ {\large and check tightenings }};
        \node (out1) [decision, below of=pro2a, yshift=-1cm , fill=Gray!30] {$max|\Omega_{V i}^{m} - \Omega_{V i}^{m-1}| \leq \eta_{V}^{\Omega} $ \\
        \& \\
         $max|\Omega_{i_{br}}^{m} - \Omega_{i_{br}}^{m-1}| \leq \eta_{I}^{\Omega} $};
        \node (stop) [startstop, below of=out1, yshift=-1cm, opacity=0.9] {{\large Stop}};
        
        \draw [arrow,fill=green] (start) -- (in1);
        \draw [arrow] (in1) -- (pro1);
        \draw [arrow] (pro1) -- (dec1);
        \draw [arrow] (dec1) -- node[anchor=east] {Yes} (pro2a);
        \draw [arrow] (dec1.west) -- +(-0.5,0) node[anchor=north] {No} |-  (in1);
        \draw [arrow] (pro2a) -- (out1);
        \draw [arrow] (out1) -- node[anchor=east] {Yes} (stop);
        \draw [arrow] (out1.west) -- +(-0.65,0) node[anchor=north] {No}  |-   (in1);
        
        
        \begin{scope}[on background layer]
                \node [note2, anchor=north west, text width=2cm] at (-4.35,-1.25) {};
        \end{scope}
        \begin{scope}[on background layer]
                \node [note, anchor=north west, text width=2cm] at (-3.5,-1.25) {};
        \end{scope}
        \node[align=center,rotate=-90,] at (3.8,-4.6) {{\LARGE Multi - period BFS-OPF}};
        
        
        \node[align=center,rotate=-90,] at (3.8,-11.25) {{\LARGE Uncertainty tightenings} };
        \node[align=center,rotate=-0,] at (0.15,-3.0) {Optimal \,  Setpoints };
        
        \node[align=center,rotate=-0,] at (0.6,-4.8) {$ V_{\textrm{bus}}^{PF} $ };
        \node[align=center,rotate=-0,] at (-2.8,-2.5) {$ V_{\textrm{bus}}^{k} $ };
        
        \node[align=center,rotate=-0,] at (-3.75,-2.5) {$ \Omega_{V i}^{m} $ };
        \node[align=center,rotate=-0,] at (-3.8,-3.1) {$ \Omega_{i_{br}}^{m} $};
        
    \end{tikzpicture}
    \caption{Proposed centralized CC-OPF scheme for the computation of the optimal DER setpoints.}
	\label{fig:SolAlgo}
\end{figure}
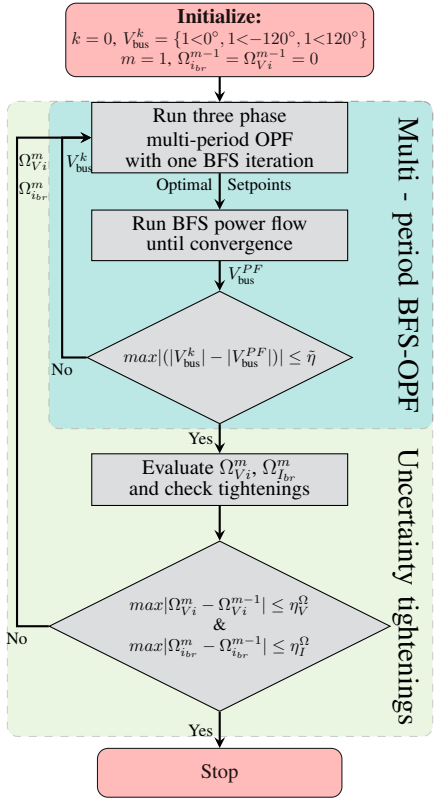

In this section, we summarize the proposed solution method for the centralized CC-OPF scheme, sketched in Fig.~\ref{fig:SolAlgo}. First, the initialization stage sets the uncertainty margins to zero and initializes the voltage levels of the three phases to a flat voltage profile. At the core of the proposed methodology lies the formulation of the three-phase multi-period centralized CC-OPF, which is summarized as 

\begin{align}
    \min_{ \bm{u} } & \sum \limits_{t=1}^{N_\textrm{hor}} \biggl\{ \sum \limits_{z\in\{a,b,c\}} \sum \limits_{j=1}^{N_\textrm{b}}  \biggl(C_\textrm{P} \hspace{-0.05cm} \cdot \hspace{-0.05cm} P_\textrm{curt,j,z,t} \hspace{-0.05cm} + \hspace{-0.05cm} C_\textrm{Q} \hspace{-0.05cm} \cdot \hspace{-0.05cm} Q_\textrm{ctrl,j,z,t} \biggr) \hspace{-0.05cm} \nonumber  \\  
    & + \sum \limits_{z\in\{a,b,c\}} \sum \limits_{i=1}^{N_\textrm{br}} C_\textrm{P} \hspace{-0.05cm} \cdot \hspace{-0.05cm} \nonumber P_\textrm{loss,i,z,t} \biggl\} \cdot \Delta t \\
    &+ C_\textrm{H} \cdot \biggl( ||\eta_\textrm{V}||_{\infty} + ||\eta_{\textrm{I}}||_{\infty} +||\eta_{\textrm{VUF}}||_{\infty} \biggr)  
\end{align}
subject to
\begin{gather*}
	P_{\textrm{inj,j,z,t}}^{\textrm{ }}=P_{\textrm{g,j,z,t}}^{\textrm{ }} - P_{\textrm{lflex,j,z,t}}^{\textrm{ }} - (P_{\textrm{B,j,z,t}}^{\textrm{ch}} - P_{\textrm{B,j,z,t}}^{\textrm{dis}})\\
    Q_{\textrm{inj,j,z,t}}^{\textrm{ }}=Q_{\textrm{g,j,z,t}}^{\textrm{ }} - P_{\textrm{lflex,j,z,t}}^{\textrm{ }} \cdot \textrm{tan}(\phi_{\textrm{load}}) + Q_{\textrm{B,j,z,t}}\\  
    I_\textrm{inj,j,z,t}= \left (\frac{(P_{\textrm{inj,j,z,t}}^{\textrm{ }} + jQ_{\textrm{inj,j,z,t}}^{\textrm{ }})^{*}}{\bar{V}_{\textrm{j,z,t}}^{*}}\right)   \\
    \bm{I}_\textrm{br,t}=\bm{BIBC} \cdot \bm{I}_\textrm{inj,t}  \\ 
    \Delta{\bm{V}}_{\textrm{t}}=\bm{BCBV} \cdot \bm{I}_\textrm{br,t} \\
    \bm{V}_{{\textrm{t}}}= \bm{V}_\textrm{slack} - \Delta V_\textrm{tap} \cdot \rho_{\textrm{t}} + \Delta{\bm{V}}_{\textrm{t}}  \\
     \rho_{min} \leq \rho_{\textrm{t}} \leq \rho_{max}   \\
     |V_{\textrm{-,j,t}}|  \leq VUF_{\textrm{MAX}} + \eta_{\textrm{VUF,j,t}},  \qquad  \eta_{\textrm{VUF,j,t}}  \geq 0\\
     P_{\textrm{g,j,z,t}}^{\textrm{min}} \leq P_{\textrm{g,j,z,t}}^{\textrm{ }} \leq P_{\textrm{g,j,z,t}}^{\textrm{max}}\\  Q_{\textrm{g,j,z,t}}^{\textrm{min}} \leq  Q_{\textrm{g,j,z,t}}^{\textrm{ }} \leq Q_{\textrm{g,j,z,t}}^{\textrm{max}}\\
    P_{\textrm{lflex,j,z,t}}^{\textrm{ }} = P_{\textrm{l,j,z,t}}^{\textrm{ }} + n_{\textrm{j,z,t}} \cdot P_{\textrm{shift,j,z}}, \qquad
   \sum \limits_{t=1}^{N_{hor}} n_{\textrm{j,z,t}}=0\\
    SoC_{\textrm{min}}^{\textrm{bat}} \cdot E_{\textrm{cap,j,z}}^{\textrm{bat}} \leq  E_{\textrm{j,z,t}}^{\textrm{bat}} \leq SoC_{\textrm{max}}^{\textrm{bat}} \cdot E_{\textrm{cap,j,z}}^{\textrm{bat}}\\  
    E_{\textrm{j,z,1}}^{\textrm{bat}} = E_{\textrm{start}}\\
    E_{\textrm{j,z,t}}^{\textrm{bat}}  = E_{\textrm{j,z,t-1}}^{\textrm{bat}} + (\eta_{\textrm{bat}} \cdot  P_{\textrm{B,j,z,t}}^{\textrm{ch}} - \frac{P_{\textrm{B,j,z,t}}^{\textrm{dis}}}{\eta_{\textrm{bat}}}) \cdot \Delta t \\
    0 \leq P_{\textrm{B,j,z,t}}^{\textrm{ch}} \leq P_{\textrm{max}}^{\textrm{bat}} ,  \quad  0 \leq P_{\textrm{B,j,z,t}}^{\textrm{dis}} \leq P_{\textrm{max}}^{\textrm{bat}} \\
    P_{\textrm{B,j,z,t}}^{\textrm{ch}} +   P_{\textrm{B,j,z,t}}^{\textrm{dis}} \leq  \textrm{max}(P_{\textrm{B,j,z,t}}^{\textrm{ch}}, P_{\textrm{B,j,z,t}}^{\textrm{dis}}) \\
    Q_{\textrm{B,j,z,t}}^{2} \leq  (S_{\textrm{max}}^{\textrm{bat}})^{2} - \textrm{max}((P_{\textrm{B,j,z,t}}^{\textrm{ch}})^{2}, (P_{\textrm{B,j,z,t}}^{\textrm{dis}})^{2})\\
    \begin{cases}
               |\mathcal{R} V_{\textrm{j,z,t}} |  \leq V_{\textrm{max}} - \Omega_\textrm{V j,z,t}^\textrm{upper} + \eta_{\textrm{V,j,z,t}} \\
               {Re}\left\{\mathcal{R} V_{\textrm{j,z,t}}    \right\} \geq V_\textrm{min} + \Omega_\textrm{V j,z,t}^\textrm{lower} -\eta_{\textrm{V,j,z,t}}
    \end{cases} \\
    |I_\textrm{br,i,z,t}|  \leq I_{\textrm{i,z,max}} - \Omega_{I_\textrm{br,i}} + \eta_{\textrm{I,i,z,t}}.
\end{gather*}

The BFS-OPF block calculates the optimal DER setpoints based on a single sweep of the BFS algorithm. Thus, the single iteration of the BFS equations replaces the non-convex, exact AC power flow equations with a linearized version. After we obtain the OPF setpoints, we run an \emph{exact power flow algorithm} using the obtained control settings to project the solution to the AC feasible space. The BFS-OPF block is then performed again using the updated voltages from the exact BFS power flow. These inner iterations are carried out until convergence. After the multi-period BFS-OPF has converged, we account for uncertainties in the outer loop as described in Section~\ref{uncertain}. The uncertainty margins are evaluated using the Monte Carlo approach, i.e. running AC power flows with samples of the uncertain PV injections. The iteration index of the OPF loop is denoted by $k$ and the iteration of the uncertainty loop by $m$. The iterative procedure continues until all parts of the algorithm have reached convergence.

The convergence characteristics of the proposed method are analyzed in~\cite{StavrosIREP,Roaldirep}. In~\cite{StavrosIREP}, we show that the algorithm works well for practical cases where the OPF solution does not change significantly, i.e. does not show sudden changes from iteration to iteration. In this case, the tightenings do not change much, and convergence is reached after a few iterations. However, there might be cases where the algorithm does not converge. As explained in~\cite{Roaldirep}, subsequent iterates in the algorithm might cycle between repeated points that have large differences in the associated tightenings. In that paper, the authors followed a “cut-and-branch” approach to interrupt the cycling and enforce convergence. In our case, when we are faced with non-convergent cases, the algorithm uses an ‘acceleration factor’ changing the solution less aggressively to avoid oscillations between repeated points, but at the cost of increasing the needed number of iterations.
\vspace{-0.2cm}
\section{Optimized local control design}\label{decentralized}
In this section, we describe the core idea of the paper as summarized in Fig.~\ref{fig:overview1}. First, we discuss what kind of data are needed to perform the offline CC-OPF, which is explained in detail in Section~\ref{detOPF}. Then, after deriving the optimal DER setpoints, we explain how we design the individual local controls for each DER using various ML techniques.

\begin{algorithm}[b]
\textbf{Input:} Optimal DG setpoints \\
\textbf{Output:} Optimized local characteristic curve
\begin{algorithmic}[1]
\State Set $n_{s}$, initialize the break-points $s$, $i=1$, $RSS_{\textrm{0}}=1000$ 
\State \textbf{Iterate:}
$$RSS_{\textrm{i}} := \min\limits_{\tilde{x}_{0}, \bm{\beta, s,\gamma}}\sum \limits_{t}^{T} P_{\textrm{g,t}}^{\textrm{ }}\cdot (x_{t}^{} - \tilde{x}_{t}^{ })^{2} + \sum \limits_{k=1}^{n_{s}} \gamma_{k}^{2}$$ \textrm{subject to} \newline
$\tilde{x}_{t}^{ } \hspace{-0.051cm} = \tilde{x}_{0}+\hspace{-0.051cm} \beta_{0} \cdot v_{t} + \sum \limits_{k=1}^{n_{s}} \beta_{k} \cdot (v_{t} \hspace{-0.051cm} - s_{k}^{i}) \cdot I(v_{t}\hspace{-0.051cm} > \hspace{-0.1cm} s_{k}^{i}) + \newline
-\sum \limits_{k=1}^{n_{s}} \gamma_{k} \cdot I(v_{t}\hspace{-0.051cm} > \hspace{-0.1cm} s_{k}^{i})$ \newline
 Monotonicity and slope constraints 
\State Update $s_{k}^{i+1}=\frac{\gamma_{k}}{\beta_{k}}+s_{k}^{i}$ and iteration index $i=i+1$
\State \textbf{Until:} $|RSS_{\textrm{i}} - RSS_{\textrm{i-1}}| < 0.0001$ 
\State Post-process the derived characteristic curves to be complete for all voltage values.
\end{algorithmic}
\textbf{Return:} Break-points $s$ and slope factors $\beta$ for each DG as $\{P,Q\}_{DG}\hspace{-0.1cm}=\hspace{-0.1cm}f(V)$
\caption{Local DG control design ($x\in\{p,q\}$)} \label{alg_curves}
\end{algorithm}

\subsection{Optimal DER setpoint data generation}
The first step is to generate the optimal DER setpoint data that will be used for the training of the local controls. To do this, several operating scenarios are selected from seasonal historical data. Then, the CC-OPF of Section~\ref{OPFall} is used to compute \emph{off-line} the optimal DER setpoints. The selection of the scenarios is critical, since they will form the basis for the training of the local schemes.

The DSO does not know the exact generation of all PVs in the LV system in the operational stage (as this would require detailed monitoring of all PVs). However, the DSO is aware of the installed DG capacity and the PV generation can be estimated with some uncertainty, using historical expected PV injection data and the installed capacity. These estimates are used in the CC-OPF solution. 

\review{The proposed method can be used with different seasonal data to account for seasonalities in terms of the DG injections and load. By changing the local control schemes based on the season, e.g. using the actual date, or when the topology of the DN changes, one can easily derive a behavior close to the optimal during the whole year.}

\subsection{Derivation of DG local controls}\label{deri}
For the DGs, we derive optimized local controls for Active Power Curtailment (APC) and Reactive Power Control (RPC). These controls take the form of simple, piece-wise linear characteristic curves (such as in~\cite{KaragiannopoulosIET}), much like the local control schemes used today in industry. Unlike the current industry standards, these characteristics might have an arbitrarily large number of piece-wise linear segments and are optimized for each individual DG and DN.

Defining the location of the break-points and the slope coefficients is a non-linear and non-differentiable problem. Thus, we employ the method in Algorithm~\ref{alg_curves} coming from~\cite{Muggeo2003} that iteratively refines the location of the break-points while solving a constraint residual sum-of-squares (RSS) optimization problem for the slope coefficients.

First, we define the number of break-points $n_{s}$ and initialize them. Then, we use the iterative steps $2-4$, where we solve the residual sum of squares problem using the active power injections as weights in the objective function, fitting the linear equivalent estimation model taking into account monotonicity and slope constraints. As inputs, we use the voltage $v_{t}$ for each sample $t$, $\forall t=1,\ldots,T$. Then, we fit the linear model based on the known breakpoints $s_{k}^{i}$, $\forall k=1,\ldots,n_{s}$ at the current iteration $i$, the left slope $\beta_{0}$ and difference-in-slopes $\beta_{k}$. The indicator function $I(\cdot)$ becomes one when the inside statement is true. Finally, $\tilde{x}_{0}$ is the model intercept and $\gamma$ a parameter which updates the location of the breakpoints towards the optimal one.

Omitting the indices for clarity, the key idea is to substitute the non-linear function $\beta \cdot (v \hspace{-0.051cm} - s) \cdot I(v\hspace{-0.051cm} > \hspace{-0.1cm} s)$ where both the difference-in-slope and the break-points are unknown, with its Taylor expansion using fixed break-points at each iteration $i$
$\beta \cdot (v \hspace{-0.051cm} - s^{i}) \cdot I(v\hspace{-0.051cm} > \hspace{-0.1cm} s^{i}) - \gamma\cdot I(v\hspace{-0.051cm} > \hspace{-0.1cm} s^{i})$~\cite{Muggeo2003}.

The same method is used both for the APC and RPC curves, using respectively the PV optimal active ($p$) and reactive ($q$) setpoints from the CC-OPF.

\begin{algorithm}[b]
\textbf{Input:} Optimal BESS setpoints\\
\textbf{Output:} SVM model for the real-time BESS response
\begin{algorithmic}[1]
\State For each BESS unit form $\Phi=[V,P_{\textrm{load}},Q_{\textrm{load}},P_{\textrm{g}}]$ and assume a function $f(\Phi)=\bigl\langle w,\Phi \bigr\rangle +b$.
\State Apply the linear, polynomial and radial-basis function kernels to $\Phi$. \newline
Solve: $\min\limits_{w,b,\xi } \frac{1}{2} w^T w+C\sum\limits_{t=1}^{T}(\xi +\xi^{*})$
\newline \textrm{subject to} \newline
$x - \bigl\langle w,\Phi_i\bigr\rangle - b\leq\epsilon+\xi,  \forall \,\,\,  (\Phi_{i},x) \newline
\bigl\langle w,\Phi_i\bigr\rangle + b - x \leq \epsilon+\xi^{*},  \forall \,\,\, (\Phi_{i}, x)$ 
\State Identify the kernel with the lowest out-of-sample error
\end{algorithmic}
\textbf{Return:} $\{P,Q\}_{B}\hspace{-0.1cm}=\hspace{-0.1cm}f(V,P_{\textrm{load}},Q_{\textrm{load}},P_{\textrm{g}})$
\caption{Local BESS control design ($x\in\{p,q\}$)} \label{Bess_local}
\end{algorithm}

\subsection{Local control of Battery Energy Storage Systems}\label{bessdecentralized}
Due to the more complex behavior of BESSs, e.g. inter-temporal constraints, we chose an SVM regression model to approximate the optimal setpoints of active and reactive power for the BESSs. An SVM regression model calculates a function, which deviates from the training data by a value no greater than a predetermined margin ($\epsilon$ in Algorithm~\ref{Bess_local}), and at the same time is as flat as possible. SVM models are very powerful because they can also model nonlinear functions (or decision boundaries). This is achieved, by mapping the training set from the input space into higher dimensional spaces, called feature spaces, by performing a non-linear transformation using suitably chosen basis functions (kernels). Then, they solve the linear model in the new space problem, which describes a nonlinear behavior in the original (input) space~\cite{hastie_09_elements-of.statistical-learning}. The procedure of training the SVM controllers follows~\cite{FedericaGM} and is summarized in Algorithm~\ref{Bess_local}. We use as features ($\Phi$) the local active and reactive power demand ($P_{\textrm{load}}, Q_{\textrm{load}}$), the active power injection of the PV at the same node ($P_{\textrm{g}}$) and the local voltage measurement ($V$). We then use these local features, in their actual or a higher dimensional space through Kernels, to create a model that mimics the optimal response by the CC-OPF setpoints.

In order to derive the best SVM model, we test three different Kernels: the linear ($\bigl\langle \Phi,\Phi^T\bigr\rangle$ in which case $C$ in the objective function is a free parameter), the polynomial ($(\gamma \bigl\langle \Phi,\Phi^T\bigr\rangle +r)^d$ where $C$ and the polynomial order $d$ are free parameters) and the Radial-Basis Function (RBF) Kernel ($e^{(-\gamma |\Phi - \Phi^T|)^2}$ where $C$ and the kernel scale $\gamma$ are free parameters). Assuming a regression function $f(\Phi)=\bigl\langle w,\Phi \bigr\rangle +b$, we solve the  convex optimization problem shown in Step~$2$, for all these Kernels in order to identify the most suitable one. The constant $C$, also called box constraint, takes positive values and penalizes the observations that lie outside the region defined by $\epsilon$, helping to prevent overfitting (regularization). The value of $C$ assesses the trade-off between the flatness of the regression function and the amount up to which deviations larger than $\epsilon$ are tolerated. 
 Finally, we keep the model with the kernel resulting in the lowest overall out-of-sample error through a 5-fold cross validation process.

\subsection{Local control of Controllable Loads}\label{CLdecentralized}
For the controllable loads, we use an SVM model as a classifier, where we define three classes $y_c \in \{-1,0,1\}$ for the `load decrease', `no shifting' and `load increase' cases, respectively. As features we use $\Phi=[V,P_{\textrm{g}}]$, where $V$ is the local voltage measurement and $P_{\textrm{g}}$ the active power injection of the PV at the same node. Intuitively, high PV injections increase local voltages triggering a load increase action.  
The optimization problem is similar to the BESS case, with $G(\Phi) = sign(f(\Phi))$ being the classifier. 

SVMs are able to deal with datasets with imbalanced class observations. This can be done by assigning different values for the box constraint (constant $C$) for `positive’ and `negative’ classes, i.e. changing the misclassification penalty for each class. This is equivalent to changing the class observation frequencies, i.e. oversampling the minority class. For example, if $C_{pos}= 2 \cdot C_{neg}$ this is in principle equivalent to training a standard SVM with $C=C_{neg}$ after considering the positive training samples twice. Such an approach implementing the so-called class-weighted SVM has been introduced in many references such as~\cite{Osuna1997}. In our model, the classes of `load increase' and `load decrease' are balanced by design, due to constraints \eqref{eq:CL}, that impose preservation of the daily load. Thus, since only the class ‘no load shifting’ can comprise different amount of observation samples, we considered a different weight in this classification class.

\section{Case Study - Results}\label{case}

\begin{figure}[]
    \begin{centering}
	\includegraphics[width=0.9\columnwidth]{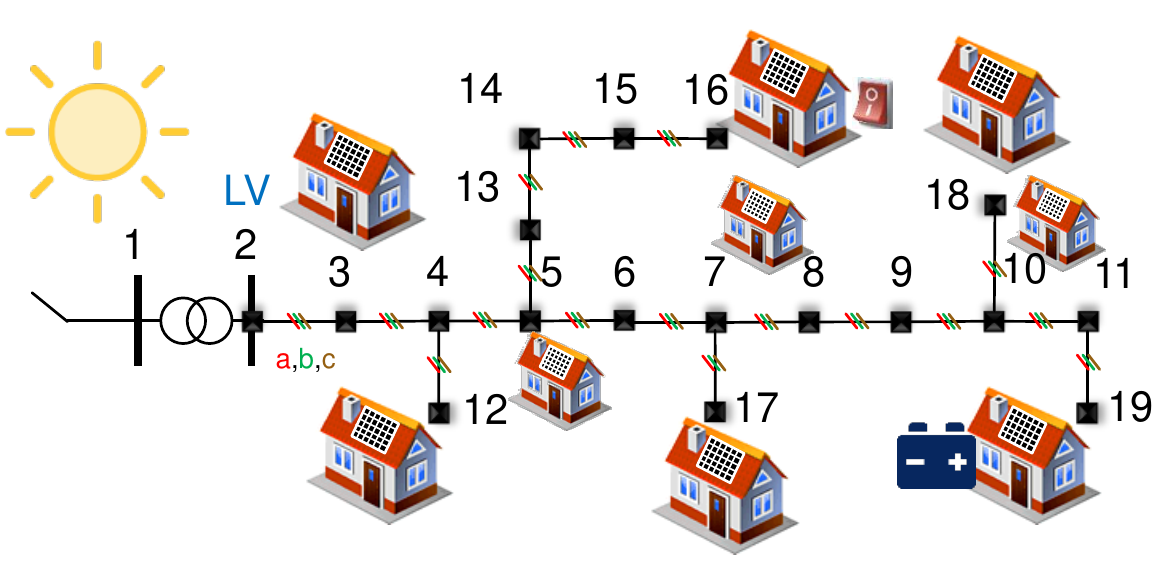}
	 \caption{Typical residential European LV grid~\cite{Strunz2014}.}	\label{fig:cigre_test_system} 
	\end{centering}
	\vspace{-0.2cm}
\end{figure}

To analyze the performance of the proposed control design algorithm, we use a typical European radial LV grid~\cite{Strunz2014}, sketched in Fig.~\ref{fig:cigre_test_system}. The neutral is assumed to be earthed in several points, and due to the short lengths of cables the capacitance is neglected. \review{The pole grounding impedance is $40~\Omega$ corresponding to distributed neutral earthing, and the transformer grounding impedance $3~\Omega$. Following~\cite{Strunz2014}, the effect of the ground return path is considered in the primitive impedance matrices. The interested reader is referred to~\cite{Kersting2002} for modeling details.}

The load and PV panels are distributed to the three phases unevenly, in order to simulate unbalanced conditions. More specifically, the total load taken from~\cite{Strunz2014} is shared $25\%$-$60\%$-$15\%$ among the three phases. The installed PV capacity, is set to $S_{\textrm{rated}}^{\textrm{PV}}=150\%$ of the total maximum load of the entire feeder to the PV nodes~$=[3, 5, 7, 10, 12, 16, 17, 18, 19]$, and is shared $25\%$-$25\%$-$50\%$ among the three phases. 

Furthermore, a BESS is located on node $19$ of phase $C$ with capacity $\frac{1}{2}S_{\textrm{rated}}^{\textrm{PV}}$~kWh, where $S_{\textrm{rated}}^{\textrm{PV}}$ is the rated power of the PV unit at that particular node. A flexible load of $5$ kW connected to phase C of Node $16$, whose total daily energy consumption needs to be constant. Please note that we assume single-phase connections for both the loads and the PV panels.

For comparison, we perform three different investigations for the operation of the system:
\begin{itemize}
    \item Method 0: The DGs are operating according to the German grid-code rules~\cite{VDE}, and no other DERs are allowed to be controlled by the DSOs. This corresponds to the current practice in industry.
    \item Method 1: All DERs are controlled based on the OPF-based algorithm described in Section~\ref{OPFall} assuming perfect communication and monitoring infrastructure. As this serves as the benchmark of the best achievable performance, we consider perfect measurements and predictions for the whole time horizon without any uncertainty.
    \item Method 2: All DERs are operating according to the individual controls derived in Section~\ref{decentralized}.
\end{itemize}

The implementation was done in MATLAB. For the centralized OPF-based control, YALMIP~\cite{Lofberg2004} was used as the modeling layer and Gurobi~\cite{gurobi} as the solver. The results were obtained on an Intel Core i7-2600 CPU and 16 GB of RAM.
\subsection{Derived local control}\label{loc_control_results}
To derive the local control schemes of all DERs, we use a 30-day summer dataset with forecasts of the PV production with 1-hour time resolution. Thus, for all cases, the training set comprises $30*24=7200$ samples. Then, the algorithm described in Section~\ref{OPFall} is used to generate the optimal DER setpoint data. The operational costs are assumed to be $C_{\textrm{P}} = 0.1 \frac{\textrm{CHF}}{\textrm{kWh}}$ and $C_{\textrm{Q}} = 0.01 \cdot C_{\textrm{P}}$. The BESS, CL, and OLTC costs are considered in the planning stage~\cite{stavrosPowertech} and their use does not incur any operational cost to the DSO. Finally, $C_{\textrm{H}}=1000 \cdot C_{\textrm{P}}$ is used to avoid infeasible solutions. For the CC-OPF, we use forecast error distributions from~\cite{Meteoswiss}, and draw $1000$ samples from the 9-hour ahead forecast error distribution of the summer power profiles similar to~\cite{StavrosIREP}. We assume a perfect spatial correlation, implying that all PVs follow the same distribution. An acceptable violation probability of $\epsilon=5\%$ is used.
Then, from the generated optimal DER setpoint data, we derive the local controls as described in Section~\ref{decentralized}.
\begin{figure}[]
    \begin{centering}
    \vspace{-0.3cm}
	\includegraphics[width=0.99\columnwidth]{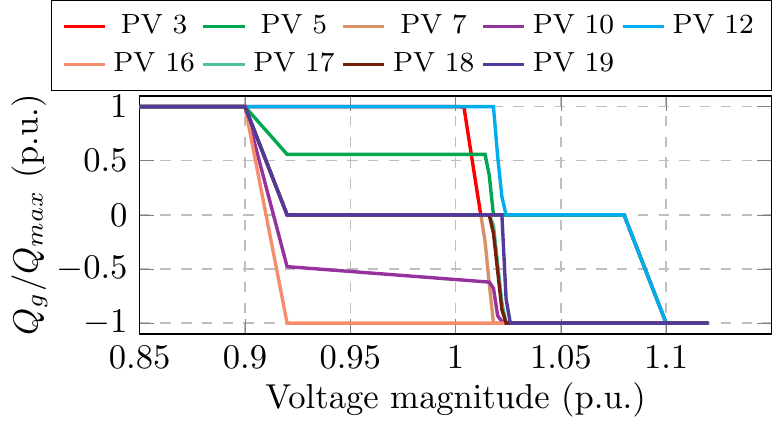}
	 \caption{Individual local characteristic curves for reactive power control of the PV units at phase C.}	\label{fig:Qcurves} \vspace{-0.3cm}
	\end{centering} 
\end{figure}
\begin{figure}[b]
    \begin{centering}
    \vspace{-0.4cm}
	\includegraphics[width=0.99\columnwidth]{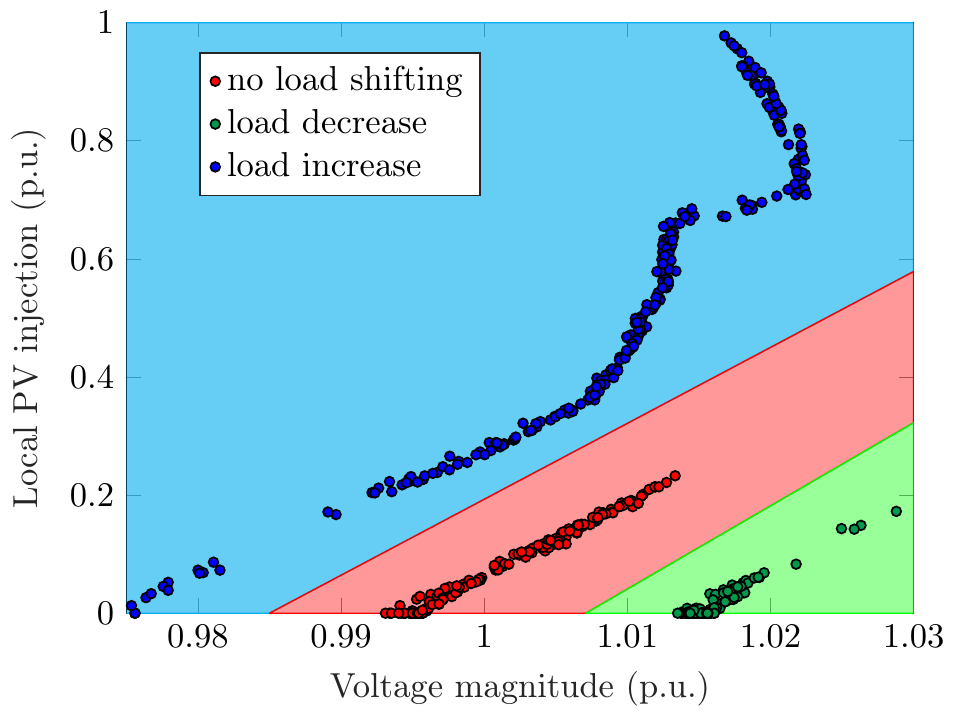}
	\caption{Classification regions defining the real - time response of the controllable loads with the x-axis indicating the local voltage magnitude and the y-axis the local PV injection.}	\label{fig:CLclasses}
	\end{centering} 
\end{figure}
Figure~\ref{fig:Qcurves} shows the individual local characteristic curves derived with Algorithm~\ref{alg_curves} for the RPC of the PV units in phase C. It can be seen that the units closer to the secondary of the substation, i.e. $3$, $5$ and $12$ show a capacitive behavior optimizing the losses, while the ones facing overvoltage problems at the end of the feeder, e.g. $16$ and $19$, show an inductive behavior at smaller voltages than the maximum of $1.04$ p.u. 

To obtain BESS and controllable load models that behave well on unseen data, we perform cross-validation (or out-of-sample testing), which is a re-sampling procedure to test the model’s performance on new data. This procedure helps identifying overfitting or selection bias issues and to provide intuition on how the model generalizes to an independent dataset. In both SVM models, we followed a 5-fold cross-validation procedure, partitioning the sample data into $5$ sets of $144$ samples. We train the SVM models using $4$ folds, and use the remaining to measure the performance.  Finally, after combining (averaging) the results of multiple rounds of cross-validation, we derive a more accurate estimate of model prediction performance.
For the BESS models, derived with Algorithm~\ref{Bess_local}, the RBF kernel functions resulted in the best behavior in terms of out-of-sample validation procedure with the following parameters:
constant $C=236.78$, $\epsilon=0.0025$, $\gamma=1.12$, showing an overall $RMSE$ of $0.158$.

Finally, for the controllable loads, the method detailed in Section~\ref{CLdecentralized} gives a classifier with overall accuracy of $100\%$ since the data are perfectly linearly separable. Figure~\ref{fig:CLclasses} shows the decision boundaries that define the three classes in the space of the two features. As can be observed, for PV injections higher than $0.3$ p.u. the load is increased to reduce the local voltage in combination with the other available measures. In low PV injections, e.g. during evening hours, the load is decreased to maintain the total daily demand constant, or is not shifted.

\vspace{-0.3cm}
\subsection{Results}\label{evaluation_results}
\begin{figure}[]
    \begin{centering}
	\includegraphics[trim=3cm 9.3cm 3.5cm 9.1cm, clip=true, width=0.9\columnwidth]{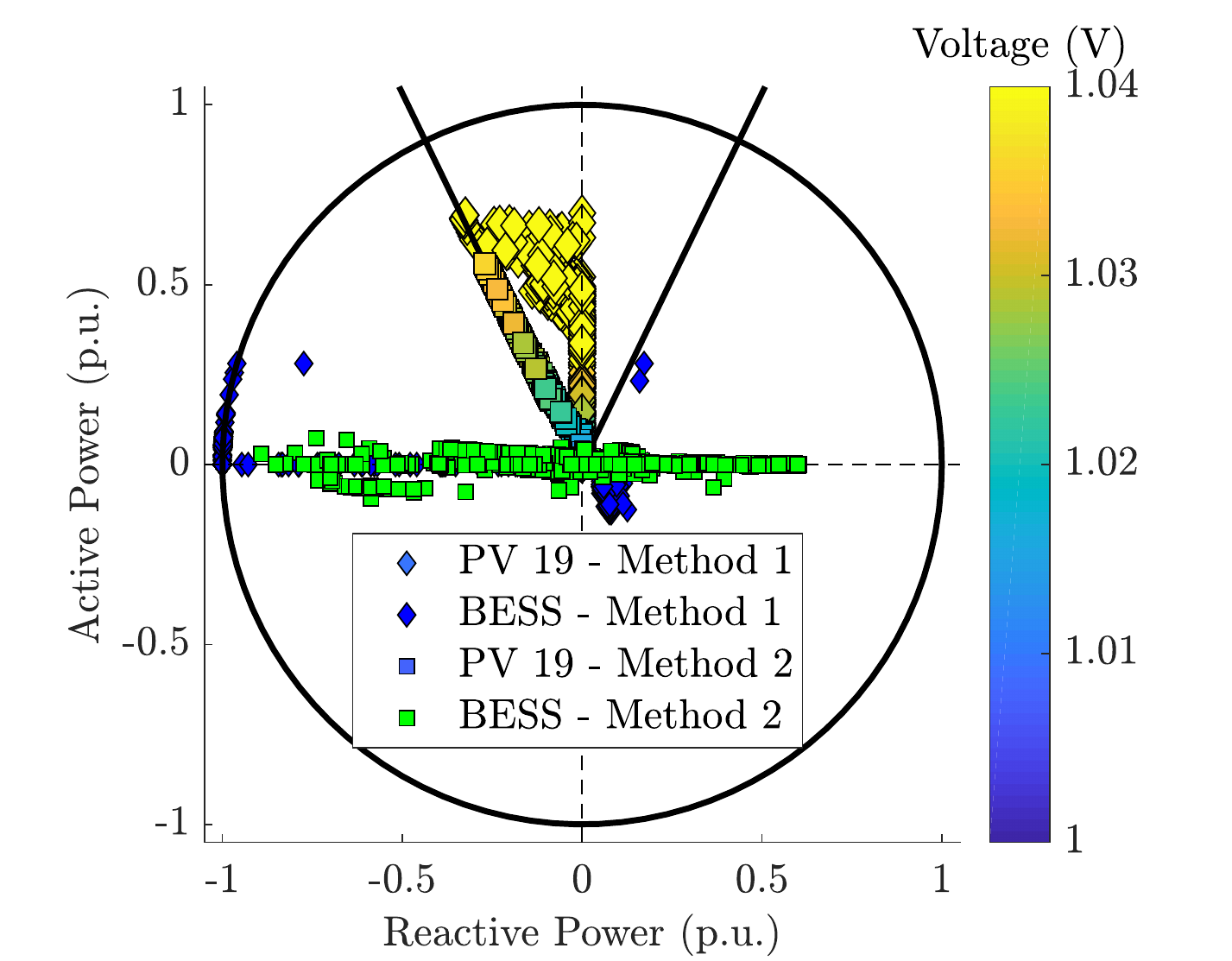}
	 \caption{Real-time control schemes of the BESS and PV unit of Node 19 phase C.}	\label{fig:BESS} 
	\end{centering}
\end{figure}

Table~\ref{SumResults} summarizes the results from applying the three methods in real-time operation for a test period of one month. Method~$1$ corresponds to the benchmark as it satisfies all security constraints and minimizes the objective function. Method~$0$ (standard industry practice) results in higher losses than the OPF-based approach, due to increased needs for reactive power by the PV units, without solving the overvoltage, overload, or balancing issues. Finally, Method~$2$ mitigates adequately the overvoltages and overloads to values acceptable by grid codes, while being capable of mimicking the OPF-based control without the need of communication. Moreover, it significantly improves the balancing problems with only small violations during $5$ hours in the month, which is also acceptable as defined by the grid codes.
\begin{table}[t]
\centering
\caption{Summarized monthly results for all methods (only the largest observed value is listed)}
\label{SumResults}
\begin{tabular}{c|ccc}
\toprule
Method      &  0     &  1 &  2 \\ \hline
Losses (\%) & 4.60   & 4.42      & 4.45       \\
$|V|_{\textrm{max}}$ (p.u.) & 1.069    & 1.04      & 1.045      \\
$|I|_{\textrm{max}}$ $(\%)$ & 119.94    & 100       & 99.49   \\
$VUF_{\textrm{max}}$ $(\%)$ & 1.81      & 1.98      & 2.33    \\
$P_{\textrm{curt}}$ $(\%)$  & 0         & 1.08      & 2.03      \\ \toprule
\end{tabular}
\end{table}
Figure~\ref{fig:BESS} displays the real-time control behaviour of the BESS and PV unit at Node 19, phase C, operating according to Methods $1$ and $2$. 
It can be seen that the proposed local control (Method~$2$) of the BESS and PV is more conservative than the OPF-based approach (Method~$1$) where the PV unit absorbs the maximum reactive power for most voltage levels due to the overvoltage problems. This conservative behavior is due to the CC-OPF approach we employ (Section~\ref{OPFall}) to generate the data used for deriving the controls of Method~$2$. On the contrary, the OPF-based control (Method~$1$) uses the actual data, assuming full knowledge of the network, load and production values through an ideal communication system without delays. Despite this, the response using the proposed local controls mimics the OPF response in a satisfactory way.

\begin{figure}[]
    \begin{centering}
	\includegraphics[width=0.95\columnwidth]{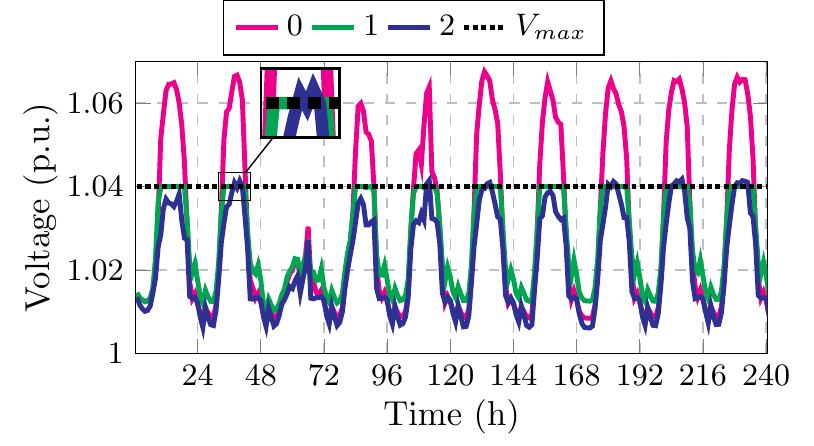}
	 \caption{Voltage magnitude evolution at phase C of Node 19.}	\label{fig:V19} 
	\end{centering}
\end{figure}

\begin{figure}[t!]
    \begin{centering}
	\includegraphics[width=0.95\columnwidth]{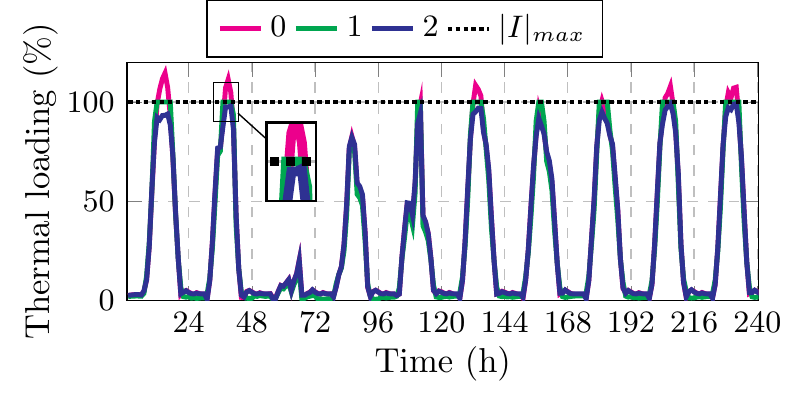}
	 \caption{Current magnitude evolution at phase C of Cable~\mbox{2--3}}	\label{fig:I23} 
	\end{centering}
\end{figure}

Finally, Fig.~\ref{fig:V19} shows the evolution of the voltage at Node~19 over the ten days. It can be seen that operating with the current regulations (Method~$0$) leads to frequent overvoltages. On the contrary, the OPF-based approach (Method~$1$) and the proposed local control satisfy the voltage security constraints. Similar observations can be made for the thermal loading, shown for Cable~\mbox{2--3} in Fig.~\ref{fig:I23}.

\section{Conclusion}  \label{Conclusion}
Future DNs will increasingly rely on the active control of DERs for the security, reliability, and optimal operation of the grid. While centralized, OPF-based controllers can provide optimal operation, they rely on expensive monitoring and communication infrastructure -- currently not available in most DNs. At the same time, the inexpensive, traditional, local controllers cannot cope with the rapidly changing environment and increased DER penetration.

In this paper, we propose a data-driven local control design methodology to derive local DER controls that can mimic the centralized controller optimal behavior, without the need for monitoring and communication infrastructure. This is based on using ML techniques to derive optimized local controls based on historical data processed through a CC-OPF. The controllers are simple to compute, understand, and implement. Yet, we have shown through the examples used that the proposed local controls can tackle security problems in an unbalanced and challenging environment while at the same time optimize its operation. Future work will focus on comparing data-driven control schemes using other ML techniques, and subsequently, on assessing risks and challenges of using such schemes at operating conditions, which were not seen in the training dataset.
\vspace{-0.4cm}
\bibliographystyle{IEEEtran}
\bibliography{bibliography}
\vspace{-0.13cm}

\vspace{-4cm}
\begin{IEEEbiography}[{\includegraphics[trim={4.5cm 14.5cm 6cm 3.8cm},width=1in,height=1.25in,clip,keepaspectratio]{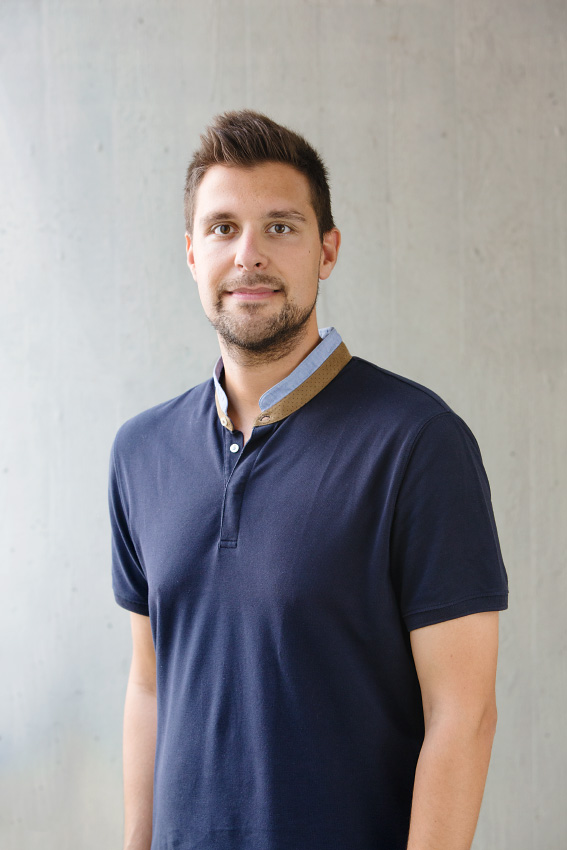}}]{Stavros Karagiannopoulos}
(S'15) was born in Thessaloniki, Greece. He received a Diploma in Electrical and Computer Engineering from the Aristotle University of Thessaloniki, Greece, in 2010, and a M.Sc. degree in Energy Science and Technology from the Swiss Federal Institute of Technology (ETH), Zurich, Switzerland, in 2013. After that, he worked at ABB Corporate Research Center in Switzerland, and since 2015, he has been pursuing the Ph.D. degree with the Power Systems Laboratory, ETH Zurich, Zurich, Switzerland. His main research focuses on planning and operation of active distribution grids.  
\end{IEEEbiography}
\begin{IEEEbiography}[{\includegraphics[trim={2.5cm 2.2cm 2.5cm .3cm}, width=1in,clip,keepaspectratio]{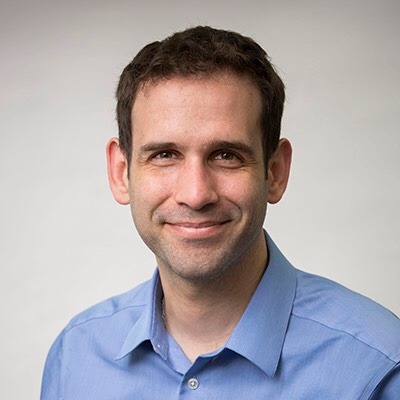}}]{Petros Aristidou}
(S'10-M'15) received a Diploma in Electrical and Computer Engineering from the National Technical University of Athens, Greece, in 2010, and a Ph.D. in Engineering Sciences from the University of Li{\`e}ge, Belgium, in 2015. He is currently a Lecturer (Assistant Professor) in Smart Energy Systems at the University of Leeds, U.K. His research interests include power system dynamics, control, and simulation.
\end{IEEEbiography}
\begin{IEEEbiography}
    [{\includegraphics[trim={5.5cm 15cm 6cm 4.25cm}, width=1in,clip,keepaspectratio]{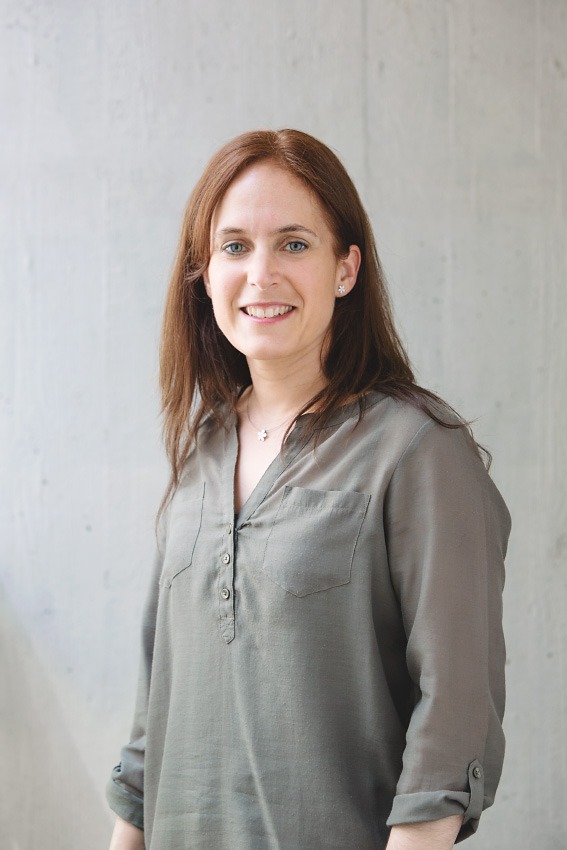}}]{Gabriela Hug}
(S'05-M'08-SM'14) was born in Baden, Switzerland. She received the M.Sc. degree in electrical engineering in 2004 and the Ph.D. degree in 2008, both from the Swiss Federal Institute of Technology, Zurich, Switzerland. After the Ph.D. degree, she worked with the Special Studies Group of Hydro One, Toronto, ON, Canada, and from 2009 to 2015, she was an Assistant Professor with Carnegie Mellon University, Pittsburgh, PA, USA. She is currently an Associate Professor with the Power Systems Laboratory, ETH Zurich, Zurich, Switzerland. Her research is dedicated to control and optimization of electric power systems.
\end{IEEEbiography}

\end{document}